%% file: Coupling_res-s_near_ess_spectrum.ltx
\documentclass{amsart}
\usepackage{amssymb}
\usepackage{mathrsfs}
\usepackage{euscript}
\usepackage{lipsum}

\usepackage{enumerate}
\input Macros.tex
\input MyMakeLit.tex

\rndef{\iff}{\Leftrightarrow}

\renewcommand{\TrD}{\mathrm{Tr}_\omega}

\begin{document}

\title[Coupling resonances near essential spectrum]{Spectral flow inside essential spectrum III:
coupling resonances near essential spectrum}

\author{Nurulla Azamov}

\address{Independent scholar, Adelaide, SA, Australia}

\email{azamovnurulla@gmail.com}
%     Secondary 47A70, 81U99.
%Primary 47A55; % Perturbation theory
%Secondary 47A11% Local spectral properties

\begin{abstract}
Given a self-adjoint operator $H_0$ and a relatively $H_0$-compact self-adjoint operator $V,$
the functions $r_j(z) = - \sigma_j^{-1}(z),$ where $\sigma_j(z)$ are eigenvalues of the compact operator $(H_0-z)^{-1}V,$ 
 bear a lot of important information about the pair $H_0$ and $V.$ These functions also admit many other characterisations, for instance, they are poles of the meromorphic function 
 $(H_0 + s V - z)^{-1}V.$  of the coupling variable $s,$ and as such they are also poles of the scattering matrix of the pair $H_0+sV$ and $H_0$ for a fixed value of energy $z.$ 
 For this reason they are called coupling resonances. 
 
 In case of rank one (and positive) perturbation $V,$ there is only one coupling resonance function, which -- in this case -- is a Herglotz function. The rank one case has been studied in depth in the literature, 
 and appears in different situations, such as Sturm-Liouville theory, random Schr\"odinger operators, harnomic and spectral analyses, etc. 
 The general case is complicated by the fact that the resonance functions are no longer single valued holomorphic functions, and potentially can have quite an erratic behaviour, typical for infinitely-valued 
 holomorphic functions. 
 
 Of very special interest are those coupling resonance functions $r_z$ which approach a real number $r_{\lambda+i0}$ from the interval $[0,1]$ as the spectral parameter $z=\lambda+iy$ approaches
 a point $\lambda$ of the essential spectrum, since they are responsible for spectral flow through $\lambda$ inside essential spectrum when $H_0$ gets deformed to $H_1 = H_0+V$ via the path $H_0 + rV, r \in [0,1].$
 
 In this paper it is shown that if the pair $H_0,$ $V$ satisfies the limiting absorption principle, then the coupling resonance functions are well-behaved near the essential spectrum in the following sense. 
 Let $I$ be an open interval inside the essential spectrum of $H_0$ and $\eps>0.$ Then there exists a compact subset~$K$ of~$I$ such that $\abs{I\setminus K} < \eps,$ and $K$ has a ``non-tangential'' neighbourhood in the upper complex half-plane, such that any coupling resonance function is either single-valued in the neighbourhood, or does not take a real value in the interval $[0,1].$
\end{abstract}
%\begin{center} { \large \sf \today \\ Draft v\,6.00} 
%\end{center}

\maketitle

%\tableofcontents

\setcounter{tocdepth}{2}
%level -1: part, 0: chapter, 1: section, etc.

%\tableofcontents

\section{Introduction}

We start by recalling definition of a coupling resonance point. For details, see \cite{AzSFIES,AzSFnRI,AzDaMN,AzDa4}.
Given a self-adjoint operator~$H_0$ and a relatively~$H_0$-compact self-adjoint operator $V,$
a \emph{coupling resonance point}, or a \emph{coupling resonance function}, $r_z,$ of the pair~$H_0$ and $V$ can be defined as a pole of the meromorphic function of the complex variable $s$ defined by 
$$
   (H_0 + sV - z)^{-1} V.
$$
Since this expression also depends on the spectral parameter~$z,$ a coupling resonance $r_z$ is a function of~$z.$

Coupling resonance points should not be mixed up with energy resonance points, see e.g.~\cite{Zw}. Both share the property of being a pole of the scattering matrix:
an energy resonance is a pole of the scattering matrix $S(z; H_0+V, H_0)$ treated as a function of the complex energy variable~$z,$ while a coupling resonance $r_z$
is a pole of the scattering matrix $S(z; H_0+s V, H_0)$ treated as a function of the complex coupling variable $s$ for a fixed value of energy~$z.$ This explains the terminology. 
One remark might be in order here. The scattering matrix $S(z; H_0+V, H_0),$ as initially defined as a function on the resolvent set of $H_0,$ cannot in fact have poles. But if it admits analytic continuation through an interval of the essential spectrum, -- which is quite a strong assumption about the operator $H_0,$ then the analytic continuation may have poles, which are then called the energy resonances. While the coupling resonances always exist no matter what the operator~$H_0$ is. Finally, we remark that there is a connection between coupling and energy resonances: if the scattering matrix admits  analytic continuation through an interval in the essential spectrum, then so do the coupling resonances and in this case their zeros are energy resonances
(while the upper and lower half-planes, where coupling resonances are initially defined, cannot have such zeros). 

Apart from the definition given above coupling resonances admit many other equivalent descriptions, see e.g. \cite{AzSFnRI}. 

The aim of this paper is to study behaviour of the coupling resonance functions near the essential spectrum
of a self-adjoint operator. Since we do not impose any conditions on the self-adjoint operator, the essential spectrum can be very pathological.
This is a major difficulty in studying resonance points near the essential spectrum.

We are specially interested in those resonance points $r_z$ which converge to a real number as $z$ approaches to a point of the essential spectrum,
since such points are responsible for the phenomenon of spectral flow inside essential spectrum, see \cite{AzSFIES}. 

The importance of coupling resonance points was implicitly well appreciated a long time ago. A lot of research was done in the case of a rank one perturbation $V$ in which case it can be assumed to be positive without loss of generality. 
For such perturbations there is only one coupling resonance point. For rank one perturbations the self-adjoint operator being perturbed can be assumed to have simple spectrum.
In this situation the resonance point is a Herglotz function which encodes all the information about the self-adjoint operator.
This theme seems to be ubiquitous in the literature as it appears in different guises, such as in Sturm-Liouville theory, in the theory of random Schr\"odinger operators,
in operator theory, in harmonic analysis, etc, see e.g. \cite{SimTrId2}, \cite{DonBk}, \cite{AiWa}.

The difficulty of the general case is that the resonance points can have branching points and possibly other types of singularities which infinitely valued holomorphic functions can possibly have. 
The main  result of this paper sorts out these issues in a domain whose boundary contains a good chunk of the essential spectrum. As it turns out, this result suffices for many applications,
some of which will appear in forthcoming papers. 

\section{Statement of the main theorem}
All Hilbert spaces in this paper are complex and separable.
Assume that $\hilb$ is a  Hilbert space, $F \colon \hilb \to \clK$ is a closed operator from $\hilb$ to another Hilbert space~$\clK,$ called a \emph{rigging operator},
and~$H_0$ is a self-adjoint operator on $\hilb,$ such that $F$ is $\abs{H_0}^{1/2}$-compact. In this setting the sandwiched resolvent 
$$
    T_z (H_0) = F R_z(H_0) F^*
$$
is a compact-operator valued holomorphic function of the spectral parameter~$z.$ We assume that~$H_0$ obeys the weak limiting absorption principle:
for almost every $\lambda \in \mbR$  the norm limit 
$$
    T_{\lambda+i0} (H_0) 
$$
exists. About the perturbation $V$ we assume that it admits a factorisation $V = F^*JF$ where $J \colon \clK \to \clK$ is a bounded self-adjoint operator. 

Further, for a point $\lambda \in \mbR$ we define $\Theta(\lambda)$ as the set of those points $z$ in $\mbC_+$ such that either $z=\lambda$
or the argument of the number $z-\lambda$ is in $[\pi/4, 3\pi/4].$ Thus, $\Theta(\lambda)$ is a kind of ``future cone'' with vertex at $\lambda.$
For a set $A\subset \mbR$ we denote by $\Theta(A)$ the union 
$$
    \Theta(A) = \bigcup_{\lambda \in A} \Theta(\lambda).
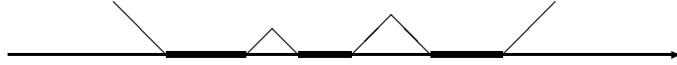
$$

\vskip 0.5 cm 
\hskip - 1.5 cm 
\begin{picture}(200,40)
\put(90,20){\vector(1,0){255}}

\put(150,20){\line(-1,1){20}}

\put(150,19){\line(1,0){30}}
\put(150,21){\line(1,0){30}}

\put(180,20){\line(1,1){10}}
\put(200,20){\line(-1,1){10}}

\put(200,19){\line(1,0){20}}
\put(200,21){\line(1,0){20}}

\put(220,20){\line(1,1){15}}
\put(250,20){\line(-1,1){15}}

\put(250,19){\line(1,0){27}}
\put(250,21){\line(1,0){27}}

\put(277,20){\line(1,1){20}}
\end{picture}

In this setting the theorem we prove here is as follows. 
\begin{thm} \label{T: main theorem} Assume the above. Let $I$ be an open interval in $\mbR$ and choose any $\eps>0.$
Then there exists a compact subset $K$ of $I$ such that $I - K$ has Lebesgue measure $< \eps$ and has the following property: 
$K$ has a neighbourhood in $\Theta(K)$ such that the restrictions of all coupling resonance points to $K$ 
either 
\begin{enumerate}
\item do not take a real value from the interval $[0,1]$ or 
\item they are single-valued continuous functions without branching points in the neighbourhood, and no they have any other kind of singularities. 
\end{enumerate}
Moreover, the number of resonance points obeying scenario (2) is finite. 
\end{thm}
We shall call the resonance points from (2) \emph{impacting}. 

The point of this theorem is that the impacting resonance points are single-valued in the chosen neighbourhood and they admit continuous extension to $K,$ which makes their analysis easy. 

Note that the functions $-r_j^{-1}(z)$ are eigenvalues of $T_z(H_0)J,$ so we may work with $r_j$ or the eigenvalue functions. 
More importantly, we may freely use properties of coupling resonances which come from the corresponding properties of isolated eigenvalues. 

\section{Sketch of proof}

\subsection{Classification of singularities of resonance functions}
A point $z_0$ in the upper half-plane we shall call \emph{singular}, if at least one coupling resonance function defined in an open set having $z_0$ in its boundary 
does not admit a single-valued analytic extension to $z_0.$ For an infinite-valued analytic function, which a collection of all resonance functions is, there can be quite a few types of singular points.
Firstly, if $z_0$ is not an isolated singular point then the number of different possible scenarios is limitless. Even if a singular point $z_0$ is isolated, there are quite a few possible scenarios:
Firstly, a resonance function can be single-valued or not in a deleted neighbourhood of $z_0,$ and if not, one can ask for the period of branching: finite or infinite. 
Secondly, the set of possible limits of $r_z$ as $z$ approaches to $z_0$ can consist of a single finite number, or of the infinity $\infty$ only, or this set can be infinite. 

As it turns out,  for coupling resonance functions there are only two possibilities given in the following theorem.
We say that a singular point $z_0,$ whether isolated or not, is \emph{absorbing} for a resonance function $r_z,$ if $z_0$ has a small enough neighbourhood~$O$ such that whenever $z$ approaches to $z_0$ along a half-interval in $O$ along which $r_z$ can be continued analytically, the limit of $r_z$ is always $\infty.$ 

\begin{thm}
Coupling resonances points can have only two types of singularities: 
\begin{enumerate}
\item isolated continuous branching points of finite period, or 
\item absorbing points, isolated or not.
\end{enumerate}
\end{thm} 
\begin{proof} Let $\gamma \in C([0,1], \rho_{ess}).$ 
If a resonance function $r(z)$ admits an analytic continuation along the half-interval $\gamma( [0,1) ),$ then the upper semicontinuity of spectrum implies that the resonance function admits continuous extension to the interval $\gamma( [0,1) ]$ in the extended complex plane. Stability of the point $r(\gamma(1))$ implies that $r(z)$ admits continuous extension to a neighbourhood of this point, possibly multivalued. 

Now, if the point $r(\gamma(1))$ is infinite, then $r(z)$ must have infinite period at this point, since otherwise $r(z)$ would take values in both half-planes in any neighbourhood of $\gamma(1),$
which a resonance function cannot.  If the point $r(\gamma(1))$ is finite, then the period of branching must be finite, as the other option contradicts compactness of the product $F T_z(H_0) F^* J.$ 
%Proof of this theorem is quite straightforward but tedious, and it is easier to prove it for oneself than to try to digest a proof written by someone else.
%The basic idea is that compactness of the product $R_z(H_0)V$ prevents all other possible scenarios. 
%For this reason we omit the proof for now.
\end{proof}

Continuous branching points of finite period exist, and it is easy to prove their existence. 
At the same time, I conjecture that absorbing points do not exist, whether isolated or not. 
Intuitively it seems unlikely that non-isolated absorbing points can exist, still 
I have been intermittently thinking about this problem for about ten years, and several times I thought I proved the conjecture only to find a gap in an argument later. It is not difficult to prove non-existence of isolated absorbing points if the sandwiched resolvent $T_z(H_0)$ is trace class. 

The following lemma shows that in any case the set of absorbing points cannot be large.
Recall that a Jordan domain is an open bounded subset of the complex plane whose boundary is homeomorphic to the circle. 

\begin{lemma} \label{L: rays lemma}   Let $G$ be a Jordan domain in the upper half-plane. Assume that $T_z(H_0)$ admits a norm-continuous extension to the closure $\bar G$
(in case the boundary $\partial G$ intersects with the real axis). Let $z_0 \in G$ and $r_z$ be a resonance function which is single-valued holomorphic in a neighbourhood of $z_0.$
Let $X$ be the set of all rays emanating from $z_0$ such that $r_z$ either does not admit single-valued analytic extension along that ray until it hits the boundary of $G$
or $r_z$ approaches $\infty$ as $z$ goes to the boundary along $\gamma.$ Then the set $X$ has Lebesgue measure zero, in an obvious sense. 
\end{lemma}
\begin{proof} This lemma follows from two classical theorems of complex and harmonic analyses: the Riemann mapping theorem which asserts that a Jordan domain is conformally equivalent to the unit disk, and that a bounded holomorphic function in the unit disk cannot non-tangentially approach to zero on a set of positive Lebesgue measure. The details are omitted. 
\end{proof}
The second theorem from the proof of this lemma is often attributed by different authors to Nevanlinna, Luzin and/or Privalov.
Also, we remark that the norm topology in this lemma can be replaced by a symmetric norm. 

\subsection{Theorem on continuous enumeration of eigenvalues}
In the proof we shall need the following theorem from \cite{AzDaTa}. 

\newcommand{\mult}{\mathrm{mult}\,}

Let $X=(X,d)$ be a based metric space with a fixed base point $x_0,$
and $\Phi$ a symmetric norm on a linear space of sequences $a = (a_1,a_2,\ldots),$ which contains $c_{00}.$
The $\Phi$-norm of a sequence $a$ we denote $\Phi(a).$
We use notation $\set{}^*$ for multisets, elements of which are counted with their multiplicities, so that for example $\set{x,x}^* \neq \set{x}^*,$ unlike the usual sets. 
We denote $\euS_\Phi(X,x_0)$
the space of all multisets $S$ in $X$ (a multiset is a set
elements $x$ of which have a multiplicity $\mult(x) = 1,2,\ldots$ which 
show how many times the element is repeated)
such that $x_0$ is an element of $S$ of infinite multiplicity and 
$$
  \Phi\brs{\set{d(x_0,x) \colon x\in S }^*} < \infty
$$
(since the norm $\Phi$ is symmetric there is no need to indicate the order in the multiset of numbers
$\set{d(x_0,x) \colon x\in S}^*$).
The support $\set{x \in S \colon \mult(x)>0}$
of a multiset from $\euS_\Phi(X,x_0)$ is countable, can accumulate only to $x_0$ 
and has no element of infinite multiplicity except $x_0.$
The space $\euS_\Phi(X,x_0)$ is endowed with a metric 
$
  \rho_\Phi(S,T) = \inf \Phi\brs{d(s_i,t_i)},
$
where the infimum is taken over all possible enumerations $S = \set{s_1, s_2, \ldots}^*$
and $T = \set{t_1, t_2, \ldots}^*$ of multisets $S$ and $T.$
In case $\Phi$ is the $p$-norm we write $\euS_p$ instead $\euS_\Phi.$ 

\begin{thm} \label{T: cont-s enumeration}   \cite[Theorem 3.1]{AzDaTa}
Let $f \colon [0,1] \to \euS_\Phi(X,x_0)$ be a continuous map.
Then there exist continuous functions $\lambda_1, \lambda_2, \ldots \colon [0,1] \to X$
such that $f(t) = \set{\lambda_1(t),\lambda_2(t), \ldots}^*$ for any $t \in [0,1].$ 
Moreover, for any neighbourhood $O$ of $x_0$ there exist only finitely many function $\lambda_j$ which take a value outside $O.$
\end{thm}

Given a symmetrically normed ideal $\clL_\Phi,(\hilb)$ the spectrum of a compact operator gives rise to a mapping 
$\clL_\Phi(\hilb) \to \euS_\Phi(\mbC, 0).$ It can be shown that this mapping is continuous, see \cite[Section 4]{AzDaTa}.
Hence, we have 
\begin{thm} \label{T: cont-s enumeration of spectrum}  
Let $f \colon [0,1] \to \clL_\Phi(\hilb)$ be a continuous map.
Then there exist continuous functions $\lambda_1, \lambda_2, \ldots \colon [0,1] \to \mbC$
such that $\sigma(f(t)) = \set{\lambda_1(t),\lambda_2(t), \ldots}^*$ for any $t \in [0,1].$ 
Moreover, for any neighbourhood $O$ of $0$ there exist only finitely many function $\lambda_j$ which take a value outside $O.$
\end{thm}
We refer to such functions as \emph{selected}.

\subsection{A version of Egorov-Severini theorem}
%\begin{thm} Assume the premise of Theorem \ref{T: main theorem}.
%Then there exists a compact subset $K_0$ of $I$ with $\abs{I-K_0} < \eps/2 $ such that  $T_z(H_0)$ admits a continuous in the norm topology extension to $\Theta(K_0).$ 
%\end{thm}

The following theorem is a modification of the classical Egorov-Severini theorem,
and it is proved by the same slightly modified argument. For reader's convenience 
we give this proof, following \cite{KF}. We recall that $\abs{A}$ denotes Lebesgue measeure of a set $A \subset \mbR.$
Also, $\Theta_\eps(\lambda)$ will denote the intersection of $\Theta(\lambda)$ with $\eps$-neighbourhood of $\lambda.$
\begin{thm} \label{T: corollary of Egorov-Severini}
Assume the premise of Theorem \ref{T: main theorem}.
Assume $T_{z}(H_0)$ takes values in a symmetrically normed ideal~${\clL_\Phi}.$
Let $\Lambda$ be a set of real numbers~$\lambda$ such that 
$T_{z}(H_0) \to T_{\lambda+i0}(H_0)$ in ${\clL_\Phi}$ as $z \to \lambda$ in $\Theta(\lambda).$
Let~$I$ be an open interval
such that $\abs{I} = \abs{I \cap \Lambda}.$
Then for any $\delta>0$ there exists a compact subset~$K$ of $I \cap \Lambda$
such that $\abs{I\setminus K} < \delta$ and 
$T_z(H_0)$ admits ${\clL_\Phi}$-continuous extension to $\Theta(K).$ 
\end{thm}
\begin{proof} We write $f(z) = T_z(H_0).$ % Let $\delta>0.$ 
For $m,n\in\mbN,$ let 
$$
  E^m_n = %\bigcap_{z \in \Theta_{1/n}(\lambda)} 
     \set{\lambda \in I \cap \Lambda \colon \forall z \in \Theta_{1/n}(\lambda) \ \ 
                     \norm{f(\lambda) - f(z)}_{\clL_\Phi} < \frac 1 m}
$$
and 
$
  E^m = \bigcup_{n=1}^\infty E^m_n.
$ 
By the premise, $\abs{I \setminus E^m} = 0.$ 
By $\sigma$-additivity of the Lebesgue measure, for any $m$ there exists $n_0(m)$ such that 
$$
  \abs{E^m \setminus E^m_{n_0(m)}} < \frac \delta{2^{m+1}}.
$$
The set $$K' := \bigcap_{m=1}^\infty E^m_{n_0(m)}$$ obeys $\abs{I \setminus K'} < \delta/2.$ 
The regularity of the Lebesgue measure allows to replace~$K'$ by a smaller compact set~$K$ such that 
$\abs{I \setminus K} < \delta.$
Further, for any $m=1,2,\ldots,$ $\lambda \in K$ and $z \in \Theta_{1/n_0(m)}(\lambda)$ 
we have 
$$
  \norm{f(\lambda)-f(z)}_{\clL_\Phi} < \frac 1 m. 
$$
Now, a small segment $[\lambda,z]$ with a slope outside of $[\pi/4,3\pi/4]$ can be replaced by 
three small segments $[\lambda,z_1],$ $[z_1,\lambda_1]$ and $[\lambda,z_2]$ with their slopes inside
$[\pi/4,3\pi/4],$ so that in addition the length of any of these three intervals is $\leq \abs{\lambda-z}.$
Hence, the $\eps/3$-trick shows that~$f$ admits a continuous extension to $\Theta(K).$
\end{proof}

%We say that a compact set $K\subset \mbR$ is an \emph{Egorov set} of 
%the pair $(H_0,F),$ if the function $T_z(H_0)$ admits a continuous extension to $\Theta(K).$ 
%Further, given $\eps>0$ we say that a compact set~$K$ is an Egorov $\eps$-set in an interval~$I$
%if~$K$ is an Egorov set with $\abs{I\setminus K} < \eps.$ 
%

\subsection{Proof of  Theorem \ref{T: main theorem}} 
By Theorem \ref{T: corollary of Egorov-Severini}, the operator function $T_z(H_0)J$ admits a continuous extension to $\Theta(K).$
By Theorem \ref{T: cont-s enumeration of spectrum}, the spectrum of the restriction of this function to the boundary of $\Theta(K)$ admits a continuous enumeration, 
$\rho_1(z), \rho_2(z), \ldots.$ Using Lemma \ref{L: rays lemma} and its proof, one can show that the set of points of $K,$ 
where at least one of the functions $\rho_j$ takes zero value, has Lebesgue measure zero.
Hence, using the method of proof of the Egorov-Severini Theorem \ref{T: corollary of Egorov-Severini}, we can assume that this null set is excluded from the compact set~$K.$
That is, none of the functions $\rho_j$ vanishes at any point of $K.$

Still, some of the functions $\rho_j$ can vanish at other points of the boundary $\Theta(K).$ But such points can be circumvented using the same Lemma 
\ref{L: rays lemma}. Thus, we can assume that none of the functions $\rho_j$ vanishes at any point of the boundary $\Theta(K).$

Since values $\rho_j(z)$ are elements of spectrum of a compact operator, this means that any $z \in \partial \Theta(K)$ has a small enough neighbourhood such that  the eigenvalue $\rho_j(z)$
can split into only a finitely many eigenvalues. Taking a compact subset of $\partial \Theta(K),$  the real part of points of which covers the interval $[-1, 2],$ we take a finite subcover of such neighbourhoods.
This shows that the set $K$ has a neighbourhood,~$G,$ in $\Theta(K),$ such that 
the selected functions $\rho_j,$ when extended from $K$ to the neighbourhood $G,$ can split only  in such a manner that the resulting functions from a finitely valued function in the neighbourhood. 

From this it is not difficult to infer that the number of branches of this multi-valued function at every point of $G$ is the same, call it $k.$
The set of branching points of this multi-valued function is discrete but it can have accumulation points in $K.$ We now show that this set of accumulation points has Lebesgue measure zero.
We choose a non-zero symmetric polynomial of $k$ variables which takes value zero if some two or more of its arguments are equal (it is not difficult to see that such a polynomial exists).
Taking composition of this symmetric polynomial with our $k$-valued holomorphic function results in a single-valued holomorphic function which vanishes at branching points (and maybe some other points too but that's ok). Applying to this function the above argument using Riemann mapping and Nevanlinna-Luzin-Privalov theorems shows that the set of branching points can accumulate only to a set of Lebesgue measure zero. Now, this null set can be cut out from the compact $K$ in the usual manner.

After that the neighbourhood $G$ of $K$ in $\Theta(K)$ will have only a finitely many branching points of our multi-valued function $\rho.$ Clearly, by making $G$ smaller, we can remove this finite set of branching points too. Proof is complete.  

%Choose a sequence of curves $\gamma_1,$ $\gamma_2,$ \ldots, each of which cut out a sort of trapezium from $\Theta(I)$
%such that they converge to the lower part of the boundary of $\Theta(K_0).$
%
%By (B), on each curve $\gamma_n$ the resonance functions admit a continuous enumeration. From these we can choose finitely many resonance functions 
%which do not take a value outside, say, 100-neighbourhood of $[0,1].$ 

%\input ../MyListOfRef
\input MyListOfRef

\end{document}

%% file: Macros.tex
\newcommand*{\abs}[1]{\left\lvert#1\right\rvert}   % absolute value of #1
\newcommand*{\set}[1]{\left\{#1\right\}}           % set of ...
\newcommand*{\brs}[1]{\left(#1\right)}             % brackets
\newcommand*{\norm}[1]{\left\Vert#1\right\Vert}    % norm of #1
\newcommand*{\hilb}{\mathcal H}                     % Hilbert space

\newcommand*{\mbC}{{\mathbb C}}
\newcommand*{\mbR}{{\mathbb R}}

\newcommand*{\mbN}{{\mathbb N}}

\newcommand*{\Tr}{\operatorname{Tr}\,}        % usual trace
    \newtheorem{thm}{Theorem}                     [section]
    \newtheorem{thm*}{Theorem}
    
    \newtheorem{lemma}[thm]{Lemma}

    \newtheorem{lemma*}{Lemma}    %for AMS \newtheorem*

    \newtheorem{rems*}{Remark}   %for AMS \newtheorem*

\newcommand{\ndef}{\newcommand*}
\def\rndef{\renewcommand}

\ndef{\myaddress}[1]{\begin{center} \it\small #1 \end{center}}

\ndef{\clA}{{\mathcal A}} \ndef{\rmA}{{\mathrm A}} \ndef{\mbA}{{\mathbb A}} \ndef{\bfA}{{\mathbf A}} \ndef{\euA}{{\EuScript A}} \ndef{\frA}{{\mathfrak A}}
\ndef{\clB}{{\mathcal B}} \ndef{\rmB}{{\mathrm B}} \ndef{\mbB}{{\mathbb B}} \ndef{\bfB}{{\mathbf B}} \ndef{\euB}{{\EuScript B}} \ndef{\frB}{{\mathfrak B}}
\ndef{\clC}{{\mathcal C}} \ndef{\rmC}{{\mathrm C}}                          \ndef{\bfC}{{\mathbf C}} \ndef{\euC}{{\EuScript C}} \ndef{\frC}{{\mathfrak C}}
\ndef{\clD}{{\mathcal D}} \ndef{\rmD}{{\mathrm D}} \ndef{\mbD}{{\mathbb D}} \ndef{\bfD}{{\mathbf D}} \ndef{\euD}{{\EuScript D}} \ndef{\frD}{{\mathfrak D}}
\ndef{\clE}{{\mathcal E}} \ndef{\rmE}{{\mathrm E}} \ndef{\mbE}{{\mathbb E}} \ndef{\bfE}{{\mathbf E}} \ndef{\euE}{{\EuScript E}} \ndef{\frE}{{\mathfrak E}}
\ndef{\clF}{{\mathcal F}} \ndef{\rmF}{{\mathrm F}} \ndef{\mbF}{{\mathbb F}} \ndef{\bfF}{{\mathbf F}} \ndef{\euF}{{\EuScript F}} \ndef{\frF}{{\mathfrak F}}
\ndef{\clG}{{\mathcal G}} \ndef{\rmG}{{\mathrm G}} \ndef{\mbG}{{\mathbb G}} \ndef{\bfG}{{\mathbf G}} \ndef{\euG}{{\EuScript G}} \ndef{\frG}{{\mathfrak G}}
\ndef{\clH}{{\mathcal H}} \ndef{\rmH}{{\mathrm H}} \ndef{\mbH}{{\mathbb H}} \ndef{\bfH}{{\mathbf H}} \ndef{\euH}{{\EuScript H}} \ndef{\frH}{{\mathfrak H}}
\ndef{\clI}{{\mathcal I}} \ndef{\rmI}{{\mathrm I}} \ndef{\mbI}{{\mathbb I}} \ndef{\bfI}{{\mathbf I}} \ndef{\euI}{{\EuScript I}} \ndef{\frI}{{\mathfrak I}}
\ndef{\clJ}{{\mathcal J}} \ndef{\rmJ}{{\mathrm J}} \ndef{\mbJ}{{\mathbb J}} \ndef{\bfJ}{{\mathbf J}} \ndef{\euJ}{{\EuScript J}} \ndef{\frJ}{{\mathfrak J}}
\ndef{\clK}{{\mathcal K}} \ndef{\rmK}{{\mathrm K}} \ndef{\mbK}{{\mathbb K}} \ndef{\bfK}{{\mathbf K}} \ndef{\euK}{{\EuScript K}} \ndef{\frK}{{\mathfrak K}}
\ndef{\clL}{{\mathcal L}} \ndef{\rmL}{{\mathrm L}} \ndef{\mbL}{{\mathbb L}} \ndef{\bfL}{{\mathbf L}} \ndef{\euL}{{\EuScript L}} \ndef{\frL}{{\mathfrak L}}
\ndef{\clM}{{\mathcal M}} \ndef{\rmM}{{\mathrm M}} \ndef{\mbM}{{\mathbb M}} \ndef{\bfM}{{\mathbf M}} \ndef{\euM}{{\EuScript M}} \ndef{\frM}{{\mathfrak M}}
\ndef{\clN}{{\mathcal N}} \ndef{\rmN}{{\mathrm N}}                          \ndef{\bfN}{{\mathbf N}} \ndef{\euN}{{\EuScript N}} \ndef{\frN}{{\mathfrak N}}
\ndef{\clO}{{\mathcal O}} \ndef{\rmO}{{\mathrm O}} \ndef{\mbO}{{\mathbb O}} \ndef{\bfO}{{\mathbf O}} \ndef{\euO}{{\EuScript O}} \ndef{\frO}{{\mathfrak O}}
\ndef{\clP}{{\mathcal P}} \ndef{\rmP}{{\mathrm P}} \ndef{\mbP}{{\mathbb P}} \ndef{\bfP}{{\mathbf P}} \ndef{\euP}{{\EuScript P}} \ndef{\frP}{{\mathfrak P}}
\ndef{\clQ}{{\mathcal Q}} \ndef{\rmQ}{{\mathrm Q}}                          \ndef{\bfQ}{{\mathbf Q}} \ndef{\euQ}{{\EuScript Q}} \ndef{\frQ}{{\mathfrak Q}}
\ndef{\clR}{{\mathcal R}} \ndef{\rmR}{{\mathrm R}}                          \ndef{\bfR}{{\mathbf R}} \ndef{\euR}{{\EuScript R}} \ndef{\frR}{{\mathfrak R}}  
\ndef{\clS}{{\mathcal S}} \ndef{\rmS}{{\mathrm S}} \ndef{\mbS}{{\mathbb S}} \ndef{\bfS}{{\mathbf S}} \ndef{\euS}{{\EuScript S}} \ndef{\frS}{{\mathfrak S}}
\ndef{\clT}{{\mathcal T}} \ndef{\rmT}{{\mathrm T}} \ndef{\mbT}{{\mathbb T}} \ndef{\bfT}{{\mathbf T}} \ndef{\euT}{{\EuScript T}} \ndef{\frT}{{\mathfrak T}}
\ndef{\clU}{{\mathcal U}} \ndef{\rmU}{{\mathrm U}} \ndef{\mbU}{{\mathbb U}} \ndef{\bfU}{{\mathbf U}} \ndef{\euU}{{\EuScript U}} \ndef{\frU}{{\mathfrak U}}
\ndef{\clV}{{\mathcal V}} \ndef{\rmV}{{\mathrm V}} \ndef{\mbV}{{\mathbb V}} \ndef{\bfV}{{\mathbf V}} \ndef{\euV}{{\EuScript V}} \ndef{\frV}{{\mathfrak V}}
\ndef{\clW}{{\mathcal W}} \ndef{\rmW}{{\mathrm W}} \ndef{\mbW}{{\mathbb W}} \ndef{\bfW}{{\mathbf W}} \ndef{\euW}{{\EuScript W}} \ndef{\frW}{{\mathfrak W}}
\ndef{\clX}{{\mathcal X}} \ndef{\rmX}{{\mathrm X}} \ndef{\mbX}{{\mathbb X}} \ndef{\bfX}{{\mathbf X}} \ndef{\euX}{{\EuScript X}} \ndef{\frX}{{\mathfrak X}}
\ndef{\clY}{{\mathcal Y}} \ndef{\rmY}{{\mathrm Y}} \ndef{\mbY}{{\mathbb Y}} \ndef{\bfY}{{\mathbf Y}} \ndef{\euY}{{\EuScript Y}} \ndef{\frY}{{\mathfrak Y}}
\ndef{\clZ}{{\mathcal Z}} \ndef{\rmZ}{{\mathrm Z}}                          \ndef{\bfZ}{{\mathbf Z}} \ndef{\euZ}{{\EuScript Z}} \ndef{\frZ}{{\mathfrak Z}}

\ndef{\tA}{{\widetilde A}} \ndef{\tcA}{{\widetilde\clA}} \ndef{\ttcA}{\widetilde{\tcA}} \ndef{\sfA}{{\textsf A}} \ndef{\ttA}{\widetilde{\tA}} \ndef{\dzA}{{A^\sharp}}
\ndef{\tB}{{\widetilde B}} \ndef{\tcB}{{\widetilde\clB}} \ndef{\ttcB}{\widetilde{\tcB}} \ndef{\sfB}{{\textsf B}} \ndef{\ttB}{\widetilde{\tB}} \ndef{\dzB}{{B^\sharp}}
\ndef{\tC}{{\widetilde C}} \ndef{\tcC}{{\widetilde\clC}} \ndef{\ttcC}{\widetilde{\tcC}} \ndef{\sfC}{{\textsf C}} \ndef{\ttC}{\widetilde{\tC}} \ndef{\dzC}{{C^\sharp}}
\ndef{\tD}{{\widetilde D}} \ndef{\tcD}{{\widetilde\clD}} \ndef{\ttcD}{\widetilde{\tcD}} \ndef{\sfD}{{\textsf D}} \ndef{\ttD}{\widetilde{\tD}} \ndef{\dzD}{{D^\sharp}}
\ndef{\tE}{{\widetilde E}} \ndef{\tcE}{{\widetilde\clE}} \ndef{\ttcE}{\widetilde{\tcE}} \ndef{\sfE}{{\textsf E}} \ndef{\ttE}{\widetilde{\tE}} \ndef{\dzE}{{E^\sharp}}
\ndef{\tF}{{\widetilde F}} \ndef{\tcF}{{\widetilde\clF}} \ndef{\ttcF}{\widetilde{\tcF}} \ndef{\sfF}{{\textsf F}} \ndef{\ttF}{\widetilde{\tF}} \ndef{\dzF}{{F^\sharp}}
\ndef{\tG}{{\widetilde G}} \ndef{\tcG}{{\widetilde\clG}} \ndef{\ttcG}{\widetilde{\tcG}} \ndef{\sfG}{{\textsf G}} \ndef{\ttG}{\widetilde{\tG}} \ndef{\dzG}{{G^\sharp}}
\ndef{\tH}{{\widetilde H}} \ndef{\tcH}{{\widetilde\clH}} \ndef{\ttcH}{\widetilde{\tcH}} \ndef{\sfH}{{\textsf H}} \ndef{\ttH}{\widetilde{\tH}} \ndef{\dzH}{{H^\sharp}}
\ndef{\tI}{{\widetilde I}} \ndef{\tcI}{{\widetilde\clI}} \ndef{\ttcI}{\widetilde{\tcI}} \ndef{\sfI}{{\textsf I}} \ndef{\ttI}{\widetilde{\tI}} \ndef{\dzI}{{I^\sharp}}
\ndef{\tJ}{{\widetilde J}} \ndef{\tcJ}{{\widetilde\clJ}} \ndef{\ttcJ}{\widetilde{\tcJ}} \ndef{\sfJ}{{\textsf J}} \ndef{\ttJ}{\widetilde{\tJ}} \ndef{\dzJ}{{J^\sharp}}
\ndef{\tK}{{\widetilde K}} \ndef{\tcK}{{\widetilde\clK}} \ndef{\ttcK}{\widetilde{\tcK}} \ndef{\sfK}{{\textsf K}} \ndef{\ttK}{\widetilde{\tK}} \ndef{\dzK}{{K^\sharp}}
\ndef{\tL}{{\widetilde L}} \ndef{\tcL}{{\widetilde\clL}} \ndef{\ttcL}{\widetilde{\tcL}} \ndef{\sfL}{{\textsf L}} \ndef{\ttL}{\widetilde{\tL}} \ndef{\dzL}{{L^\sharp}}
\ndef{\tM}{{\widetilde M}} \ndef{\tcM}{{\widetilde\clM}} \ndef{\ttcM}{\widetilde{\tcM}} \ndef{\sfM}{{\textsf M}} \ndef{\ttM}{\widetilde{\tM}} \ndef{\dzM}{{M^\sharp}}
\ndef{\tN}{{\widetilde N}} \ndef{\tcN}{{\widetilde\clN}} \ndef{\ttcN}{\widetilde{\tcN}} \ndef{\sfN}{{\textsf N}} \ndef{\ttN}{\widetilde{\tN}} \ndef{\dzN}{{N^\sharp}}
\ndef{\tO}{{\widetilde O}} \ndef{\tcO}{{\widetilde\clO}} \ndef{\ttcO}{\widetilde{\tcO}} \ndef{\sfO}{{\textsf O}} \ndef{\ttO}{\widetilde{\tO}} \ndef{\dzO}{{O^\sharp}}
\ndef{\tP}{{\widetilde P}} \ndef{\tcP}{{\widetilde\clP}} \ndef{\ttcP}{\widetilde{\tcP}} \ndef{\sfP}{{\textsf P}} \ndef{\ttP}{\widetilde{\tP}} \ndef{\dzP}{{P^\sharp}}
\ndef{\tQ}{{\widetilde Q}} \ndef{\tcQ}{{\widetilde\clQ}} \ndef{\ttcQ}{\widetilde{\tcQ}} \ndef{\sfQ}{{\textsf Q}} \ndef{\ttQ}{\widetilde{\tQ}} \ndef{\dzQ}{{Q^\sharp}}
\ndef{\tR}{{\widetilde R}} \ndef{\tcR}{{\widetilde\clR}} \ndef{\ttcR}{\widetilde{\tcR}} \ndef{\sfR}{{\textsf R}} \ndef{\ttR}{\widetilde{\tR}} \ndef{\dzR}{{R^\sharp}}
\ndef{\tS}{{\widetilde S}} \ndef{\tcS}{{\widetilde\clS}} \ndef{\ttcS}{\widetilde{\tcS}} \ndef{\sfS}{{\textsf S}} \ndef{\ttS}{\widetilde{\tS}} \ndef{\dzS}{{S^\sharp}}
\ndef{\tT}{{\widetilde T}} \ndef{\tcT}{{\widetilde\clT}} \ndef{\ttcT}{\widetilde{\tcT}} \ndef{\sfT}{{\textsf T}} \ndef{\ttT}{\widetilde{\tT}} \ndef{\dzT}{{T^\sharp}}
\ndef{\tU}{{\widetilde U}} \ndef{\tcU}{{\widetilde\clU}} \ndef{\ttcU}{\widetilde{\tcU}} \ndef{\sfU}{{\textsf U}} \ndef{\ttU}{\widetilde{\tU}} \ndef{\dzU}{{U^\sharp}}
\ndef{\tV}{{\widetilde V}} \ndef{\tcV}{{\widetilde\clV}} \ndef{\ttcV}{\widetilde{\tcV}} \ndef{\sfV}{{\textsf V}} \ndef{\ttV}{\widetilde{\tV}} \ndef{\dzV}{{V^\sharp}}
\ndef{\tW}{{\widetilde W}} \ndef{\tcW}{{\widetilde\clW}} \ndef{\ttcW}{\widetilde{\tcW}} \ndef{\sfW}{{\textsf W}} \ndef{\ttW}{\widetilde{\tW}} \ndef{\dzW}{{W^\sharp}}
\ndef{\tX}{{\widetilde X}} \ndef{\tcX}{{\widetilde\clX}} \ndef{\ttcX}{\widetilde{\tcX}} \ndef{\sfX}{{\textsf X}} \ndef{\ttX}{\widetilde{\tX}} \ndef{\dzX}{{X^\sharp}}
\ndef{\tY}{{\widetilde Y}} \ndef{\tcY}{{\widetilde\clY}} \ndef{\ttcY}{\widetilde{\tcY}} \ndef{\sfY}{{\textsf Y}} \ndef{\ttY}{\widetilde{\tY}} \ndef{\dzY}{{Y^\sharp}}
\ndef{\tZ}{{\widetilde Z}} \ndef{\tcZ}{{\widetilde\clZ}} \ndef{\ttcZ}{\widetilde{\tcZ}} \ndef{\sfZ}{{\textsf Z}} \ndef{\ttZ}{\widetilde{\tZ}} \ndef{\dzZ}{{Z^\sharp}}

\ndef{\bfa}{{\mathbf a}}
\ndef{\bfb}{{\mathbf b}}
\ndef{\bfc}{{\mathbf c}}
\ndef{\bfd}{{\mathbf d}}

\ndef{\euu}{{\EuScript u}}

  \ndef{\eps}{\varepsilon}

\let\geq\geqslant
\let\leq\leqslant

\ndef{\lims}[1]{\lim\limits_{#1}}
\ndef{\sums}[1]{\sum\limits_{#1}}
\ndef{\ints}[1]{\int_{#1}}
\ndef{\sups}[1]{\sup\limits_{#1}}
\ndef{\liminfty}[1]{\lims{#1\to\infty}}
\ndef{\suminf}[1]{\sums{#1=1}^\infty}

\ndef{\limo}[1]{\omega\mbox{-}\!\!\!\lims{#1\to\infty}}          % \omega limit
\ndef{\limL}[1]{\rmL\mbox{-}\!\!\!\lims{#1\to\infty}}            % ``L" limit
\ndef{\limLOne}[1]{\clL_1\mbox{-}\!\lims{#1}}
\ndef{\tildelimo}[1]{\tilde\omega\mbox{-}\!\!\!\lims{#1\to\infty}}
\ndef{\slim}{\mathrm{s}\mbox{-}\!\!\lim}          % strong limit
\ndef{\wlim}{\mathrm{w}\mbox{-}\!\lim}          % strong limit

\ndef{\Aut}{\operatorname{Aut}}      % group of automorphisms
\ndef{\Ch}{\operatorname{ch}}        % Chern character
\ndef{\End}{\operatorname{End}}      % group of endomorphisms
\ndef{\Hom}{\operatorname{Hom}}      % group of homomorphisms
\rndef{\ker}{\operatorname{ker}}      % kernel of an operator
\ndef{\coker}{\operatorname{coker}}      % kernel of an operator
\ndef{\im}{\operatorname{im}}        % image of an operator
\ndef{\Log}{\operatorname{Log}}      % logarithm
\ndef{\OP}{\operatorname{OP}}        % abstract PDO's
\ndef{\Op}{\operatorname{Op}}        % abstract PDO's
\ndef{\Symb}{\operatorname{Symb}}    % symbol
\ndef{\Wres}{\operatorname{Wres}}    % Wodzicki residue
\ndef{\cl}{\operatorname{cl}}        % Clifford
\ndef{\com}{\operatorname{com}}
\ndef{\const}{\operatorname{const}}  % constant
\ndef{\conv}{\operatorname{conv}}    % convex hull
\ndef{\Var}{\operatorname{Var}}
\ndef{\Cov}{\operatorname{Cov}}

\ndef{\detFK}[1]{\Delta\brs{#1}} % Fuglede-Kadison's determinant
\ndef{\detFKrel}[2]{\Delta_{#2}\brs{#1}} % Fuglede-Kadison's determinant

\ndef{\adj}{\operatorname{adj}}    % diagonal operator
\ndef{\diag}{\operatorname{diag}}    % diagonal operator
\ndef{\dist}{\operatorname{dist}}    % distance
\ndef{\dom}{\operatorname{dom}}      % domain
\ndef{\ec}{\operatorname{ec}}        % essential codimension
\ndef{\id}{\mathrm{Id}}                        % identity operator
\ndef{\ind}{\operatorname{ind}}      % index
\ndef{\mydeg}{\operatorname{deg}}    % degree of a diff. form
\ndef{\op}{\operatorname{op}}
\ndef{\rank}{\operatorname{rank}}
\ndef{\res}{\operatorname{res}}      % residue
\ndef{\ran}{\operatorname{ran}}      % range
\ndef{\sflow}{\operatorname{sf}}     % spectral flow
\ndef{\isf}{\operatorname{isf}}      % infinitesimal spectral flow
\ndef{\sign}{\operatorname{sign}}    % signum (a la C.-Ph.)
\ndef{\sgn}{\operatorname{sgn}}      % signum (a la Connes)
\ndef{\sing}{\operatorname{sing}}    % singular
\ndef{\supp}{\operatorname{supp}}    % support
\ndef{\tr}{\operatorname{tr}}        % trace
\ndef{\var}{\operatorname{var}}      % variation of measure
\ndef{\vol}{\operatorname{vol}}      % volume or volume form
\ndef{\wn}{\operatorname{wn}}        % winding number
\ndef{\wres}{\operatorname{wres}}    % Wodzicki residue
\ndef{\prng}[1]{\mathrm R_{#1}} % {[\rng {#1}]}          % projection onto the range of #1
\ndef{\pker}[1]{\mathrm N_{#1}} % {[\ker {#1}]}          % projection onto the kernel of #1
\ndef{\rprng}[2]{\mathrm R_{#1}^{#2}}           % projection onto the relative range of #1
\ndef{\rpker}[2]{\mathrm N_{#1}^{#2}}           % projection onto the relative kernel of #1
\ndef{\rsupp}[1]{\supp_r(#1)}
\ndef{\lsupp}[1]{\supp_l(#1)}
\ndef{\rslv}[1]{R_z(#1)}      % resolvent
\ndef{\HH}{H}                 % initial operator of perturbation theory
\ndef{\tHH}{\tilde \HH}       % final operator of perturbation theory
\ndef{\VV}{V}                 % perturbation operator of perturbation theory
\ndef{\Rz}{R_z}               % resolvent of the initial operator
\ndef{\tRz}{\tR_z}            % resolvent of the final operator
\ndef{\psif}[1]{#1^{[1]}} % {\psi_{#1}}  % divided difference
\ndef{\WPlus}[1]{W_{#1}(\mbR)} %\ndef{\CPlus}[1]{C^{#1+}(\mbR)}
\ndef{\bndl}{\xi}                         % vector bundle
\ndef{\bndlA}{\eta}                       % vector bundle
\ndef{\GlueMap}{\varphi}                  % glue map of a bundle
\ndef{\ChartMap}{h}                       % chart diffeomorphism map of a manifold
\ndef{\chern}{\ensuremath{\mathrm{ch}}}
\ndef{\hilba}{\clH^{(a)}}                    % absolutely continuous part of Hilbert space (wrt to a s.a. operator)
\ndef{\hilbs}{\clH^{(s)}}                    % singular part of Hilbert space (wrt to a s.a. operator)
   \ndef{\hilbasargument}{(\hilb)} %{(\hilb)}
\ndef{\LpH}[1]{\clL_{#1}\hilbasargument}       % the set of ...
\ndef{\saLpH}[1]{\clL_{sa}^{#1}\hilbasargument}       % the set of ...
\ndef{\clBH}{\clB\hilbasargument}              % the set of BOUNDED linear operators on Hilbert space
\ndef{\ubBH}{\clB_1\hilbasargument}            % the unit ball of the algebra of all BOUNDED linear operators on Hilbert space
\ndef{\clCH}{\clC\hilbasargument}              % the set of CLOSED DENSELY-DEFINED linear operators on Hilbert space
\ndef{\clKH}{\clK\hilbasargument}              % the set of COMPACT operators
\ndef{\clFH}{\clF\hilbasargument}              % the set of BOUNDED FREDHOLM operators
\ndef{\clUH}{\clU\hilbasargument}              % the set of UNITARIES on Hilbert space
\ndef{\clCFH}{{\clC\clF}\hilbasargument}       % the set of CLOSED DENSELY-DEFINED FREDHOLM OPERATORS on Hilbert space
\ndef{\saBH}{\clB_{sa}\hilbasargument}         % the set of S.-A. BOUNDED operators on Hilbert space
\ndef{\saCH}{\clC_{sa}\hilbasargument}         % the set of CLOSED DENSELY-DEFINED S.-A. operators on Hilbert space
\ndef{\saFH}{\clF_{sa}\hilbasargument}         % the set of BOUNDED FREDHOLM s.-a. operators
\ndef{\saKH}{\clK_{sa}\hilbasargument}         % the set of COMPACT S.-A. operators
\ndef{\saCFH}{\clC\clF_{sa}\hilbasargument}    % the set of CLOSED DENSELY-DEFINED S.-A. FREDHOLM operators on Hilbert space
\ndef{\clUFH}{\clU\clF\hilbasargument}         % the set of UNITARIES such that U+I is FREDHOLM
\ndef{\Uinj}{\clU_{inj}\hilbasargument}        % the set of UNITARIES such that U-I is injective
\ndef{\UFinj}{\clU\clF_{inj}\hilbasargument}   % the set of UNITARIES such that U-I is INJECTIVE and U+I is FREDHOLM

\ndef{\spproj}[2]{E^{#1}_{#2}}                      % spectral projection of #1
\ndef{\spprojb}[2]{E^{#2}_{#1}}                     % spectral projection of #1

\ndef{\LpN}[1]{\clL^{#1}(\clN,\tau)}     % noncommutative \mathcal L_p space
\ndef{\saLpN}[1]{\clL^{#1}_{sa}(\clN,\tau)} % s.-a. part of noncommutative \mathcal L_p space
\ndef{\rLpN}[1]{L^{#1}(\clN,\tau)}       % noncommutative L_p space
\ndef{\clAND}{(\clA,\clN,D)}             % spectral triple (A,N,D)
\ndef{\clBA}{{\clB(\clA)}}
\ndef{\saKN}{{\clK_{sa}(\clN,\tau)}}          % s.-a. \tau-compact operators
\ndef{\clKN}{{\clK(\clN,\tau)}}          % \tau-compact operators
\ndef{\clKtN}{{\clK(\tilde\clN,\tau)}}   % \tau-compact (maybe unbounded) operators
\ndef{\clFN}{{\clF(\clN,\tau)}}          % \tau-Fredholm operators
\ndef{\saFN}{{\clF_{sa}(\clN,\tau)}}     % self-adjoint \tau-Fredholm operators
\ndef{\clPN}{\clP(\clN)}                 % projections of N
\ndef{\clQN}{\clQ(\clN,\tau)}            % Calkin algebra N/K
\ndef{\infPN}{{\clP_\tau^\infty(\clN)}}  % infinite projections of N
\ndef{\clOF}[2]{\clF_{#1\mbox{-}#2}(\clN,\tau)}         % relatively Fredholm operators
\ndef{\oind}[2]{{\rm \tau\mbox{-}ind}_{#1\mbox{-}#2}}   % relative index
\ndef{\tind}{\tau\mbox{-}\ind}                  % semifinite index
\ndef{\DInd}{\ind_{\clD,\tau}}           % semifinite index
\ndef{\BF}{Breuer-Fredholm}              % Breuer-Fredholm
\ndef{\skewfred}[2]{$(#1\cdot #2)$ $\tau$\tire Fredholm}   % skew corner Fredholm
\ndef{\affl}{\eta}                       % affiliated
\ndef{\vNa}{von Neumann algebra}         % von Neumann algebra
\ndef{\nsf}{faithful normal semifinite } % normal semifinite faithful
\ndef{\taubrs}[1]{\tau\brackets{#1}}     % n.s.f. trace with brackets
\ndef{\sqbrs}[1]{[#1]}        % brackets
\ndef{\Sqbrs}[1]{\big[#1\big]}        % brackets
\ndef{\SqBrs}[1]{\Big[#1\Big]}        % brackets
\ndef{\domd}{\bigcap\limits_{n\ge 0} \dom\;\delta^n}         % domain of \delta^n's
\ndef{\DiffOP}{{\rm \clD}}
\ndef{\ADA}{\clA \cup [\clD,\clA]}
\ndef{\DixIdeal}[1]{\LpH{#1,\infty}}               % Dixmier ideal
\ndef{\dixideal}{\ell^{1,\infty}}                  % commutative Dixmier ideal
\ndef{\WDixIdeal}{\LpH{1,\mathrm w}}               % weak Dixmier ideal
\ndef{\DixIdealPos}[1]{\DixIdeal{#1}_+}            % positive part of Dixmier ideal
\ndef{\DixIdealN}[1]{\LpN{#1,\infty}}              % semifinite Dixmier ideal
\ndef{\DixIdealNPar}[2]{\clL^{#1,\infty}_{#2}(\clN,\tau)}    % semifinite Dixmier ideal
\ndef{\DixIdealNPos}[1]{\LpN{#1,\infty}_+}                   % positive part of semifinite Dixmier ideal
\ndef{\TrD}{\Tr_\omega}                                      % Dixmier trace
\ndef{\tauD}{{\tau_\omega}}                                  % semifinite Dixmier trace
\ndef{\ILogN}{\frac 1{\log(1+N)}}
\ndef{\DixNorm}[1]{\norm{#1}_{(1,\infty)}}                   % Dixmier norm
\ndef{\DixInt}[1]{\ints 0^t \mu_s(#1)\,ds}
\ndef{\DixIntL}[1]{\ints 0^{\lambda_{1/t}(#1)}\mu_s(#1)\,ds}
    \ndef{\SmallIdeal}{{\clL^{1, \mathrm w}}}
    \ndef{\SmallIdealMeas}{{\clL^{1, \mathrm w}_m}}
    \ndef{\DixIntII}[1]{\int_0^t \mu_s(#1)\,ds}
    \ndef{\DixIntf}[1]{\Phi_t(#1)}
    \ndef{\DixIntg}[1]{\Psi_t(#1)}

\ndef{\lpi}{\clL^{1,\pi}(\clN,\tau)}
\ndef{\strl}[1]{\stackrel \longrightarrow {#1}}
\ndef{\IIinfty}{$\mathrm{II}_\infty$\ }

\ndef{\fourier}[1]{\clF(#1)}          % Fourier transform of #1
\ndef{\HaarMeasBohrs}{\nu}            % Haar measure of the Bohr compact
\ndef{\BrownsMeas}{\mu}               % Brown's measure
\ndef{\BohrCont}[1]{\tilde{#1}}       % continuation of a function to the Bohr compact
\ndef{\APMean}{{M}}                   % mean value of a.p. function
\ndef{\CDSS}{{\clA_B}}                % Coburn-Douglas-Schaeffer-Singer's factor
\ndef{\matr}{{\rm Mat}}               % standard matrix algebra
\ndef{\seque}[1]{\ensuremath{\{#1_n\}_{n=1}^\infty}}    % sequence of numbers  a_1, a_2, ...
\ndef{\sequen}[2]{\ensuremath{\{#1_#2\}_{#2=1}^\infty}}    % sequence of numbers  a_1, a_2, ...
\ndef{\Seque}[1]{\ensuremath{\left(#1_0,#1_1,#1_2,\dots\right)}}    % sequence of numbers  a_1, a_2, ...
\ndef{\Cesaro}{H}                           % the Cesaro operator (on sequences)
\ndef{\CesaroRPlus}{M}                      % the Cesaro operator on positive semiaxis
\ndef{\Dilation}{D}                         % the dilation operator (on sequences)
\ndef{\Shift}{T}                            % the shift operator (on sequences)

\ndef{\TrNorm}[1]{\norm{#1}_1}              % trace norm of #1
\ndef{\HSNorm}[1]{\norm{#1}_2}              % Hilbert-Schmidt norm of #1
\ndef{\InftyNorm}[1]{\norm{#1}_\infty}      % uniform norm of #1
\ndef{\normQN}[1]{\norm{#1}_{\clQN}}        % Calkin norm of #1
\ndef{\clLpnorm}[2]{\norm{#2}_{\clL^{#1}}}    % $1,\infty$- trace norm of #1
\ndef{\clLnorm}[1]{\clLpnorm{1}{#1}}    % $1,\infty$- trace norm of #1

\ndef{\ccurve}{\gamma}                      % a curve in complex plane for Cauchy integral

\ndef{\Brs}[1]{\big(#1\big)}                % brackets
\ndef{\BRS}[1]{\Big(#1\Big)}                % brackets
\ndef{\scal}[2]{\left\langle #1,#2\right\rangle}               % scalar product
\ndef{\Scal}[1]{\left\langle #1\right\rangle}               % scalar product
\ndef{\precprec}{\prec\!\!\!\prec}
\ndef{\qeq}{\stackrel?=}
\ndef{\spectrum}[1]{\sigma_{#1}} %{\mathrm{Spec}(#1)}       % spectrum of an operator
\ndef{\spectruma}[1]{\sigma^{(a)}_{#1}} %{\mathrm{Spec}(#1)}       % absolutely continuous spectrum of an operator
\rndef{\emptyset}{\varnothing}                              % empty set
\ndef{\csupp}{c}                           % subscript for compactly supported functions
\ndef{\closure}[1]{\overline{#1}}
\ndef{\linspan}[1]{\mathrm{span}\,{#1}}
\ndef{\bddborel}[1]{B(#1)}                 % the space of bounded Borel functions on the measure space #1
\ndef{\charfunc}{\chi}
\ndef{\FrDer}{\euD}                        % Fr\'echet derivative
\ndef{\LieDer}[1]{\pounds_{#1}\,}          % Lie derivative
\ndef{\dds}{\left.\frac d{ds} \right|_{s = 0}}
\ndef{\ortcmp}[1]{#1^{\scriptscriptstyle \perp}}            % orthogonal complement of projection #1
\ndef{\Laplace}{\Delta}                    % Laplace operator

\ndef{\matrPQ}[3]
{
    \left(
      \begin{array}{cc}
        #1_{11} & #1_{12} \\
        #1_{21} & #1_{22}
      \end{array}
    \right)_{[#2,#3]}
}
\ndef{\margOK}{\marginpar{\bf \small OK}}

\newcounter{margcomcount}
\ndef{\margcom}[1]{\marginpar{\bf \small #1} \addtocounter{margcomcount}{1}
   \index{\indexcom{{\bf COMMENT: #1}}}}

\ndef{\mytimes}{\!\times\!}
\ndef{\sss}[1]{\subsubsection{}\label{#1}}

\rndef{\phi}{\varphi} \ndef{\OpenUnitDisk}{D}
\ndef{\RHS}{RHS}                            % right hand side
\ndef{\LHS}{LHS} %right and left hand side  % left hand side
\ndef{\ttt}{\Leftrightarrow}
\ndef{\then}{\Rightarrow}
\ndef{\tto}{\longrightarrow}
\ndef{\nno}{\nonumber\\}
\ndef{\newn}[1]{\index{#1} {\bfseries #1}}       % new notion
\ndef{\la}{\langle}
\ndef{\ra}{\rangle}
\ndef{\dbar}{{\;\bar{\phantom{o}} \!\!\!\! d}}
\ndef{\stl}[1]{\stackrel{\vbox to 0pt{\vss\hbox{$\scriptstyle #1$}}}}
\ndef{\mathcomment}[1]{{\hfill \qquad\qquad\qquad\text{by (#1)}}}        % for comments at the ends of lines with math formulas
\ndef{\mathcomm}[1]{{\hfill \qquad\qquad\qquad\qquad\qquad\text{#1}}}        % for comments at the ends of lines with math formulas
\ndef{\details}[1]{\smallskip\begin{center} {\bf Here:}
#1\end{center}\medskip} \ndef{\indexcom}[1]{ --- #1}
\ndef{\longsim}{\ \sim \ }              % for use in formulas.
\ndef{\tire}{-}              % for use in formulas.
\ndef{\intinfinf}{\int_{-\infty}^\infty}
     \ndef{\npartial}{\slash\!\!\!\partial}
     \ndef{\Heis}{\operatorname{Heis}}
     \ndef{\Solv}{\operatorname{Solv}}
     \ndef{\Spin}{\operatorname{Spin}}
     \ndef{\SO}{\operatorname{SO}}
     \ndef{\Index}{\operatorname{index}}
             \ndef{\p}{\partial}
             \ndef{\dd}{|\clD|}
             \ndef{\n}{\parallel}

%% file: MyMakeLit.tex
\let\LatexCite=\cite  % just renaming

\let\ifnumref\iffalse % use this command if you want LETTER REFERENCES like [CM]

\ndef{\ifuncited}[4]{\expandafter\ifx\csname used#4\endcsname\relax}

\ndef{\ifcited}[4]{\expandafter\ifx\csname used#4\endcsname\relax\else}
  \ndef{\papertitle}[1]{ \emph{#1}, }
  \ndef{\paperauthor}[2]{#2}  %{\ifnum#1=0$^*$\fi#2}
  \ndef{\pbbi}[9]{%
      \ifcited{#1}{#2}{#3}{#5}%
        \ifnumref%
          \bibitem{#5}\paperauthor{#1}{#6},\papertitle{#7}#8.%
        \else%
          \advance #9 by 1%
          \ifnum#9<1%
            \bibitem[#4]{#5}\paperauthor{#1}{#6}, \papertitle{#7}#8.%
          \else%
            \bibitem[#4{\the#9}]{#5}\paperauthor{#1}{#6},\papertitle{#7}#8.%
          \fi%
        \fi%
      \fi%
  }
  \ndef{\mbbi}[8]{%
     \ifcited{#1}{#2}{#3}{#5}%
        \ifnumref%
          \bibitem{#5}\paperauthor{#1}{#6},\papertitle{#7}#8.%
        \else%
          \bibitem[#4]{#5}\paperauthor{#1}{#6},\papertitle{#7}#8.%
        \fi%
     \fi%
  }
\ndef{\AddCite}[1]{%
   \ifuncited{0}{0}{0}{#1}%
     \expandafter\gdef\csname used#1\endcsname {}%
   \fi%
}

\def\ProcessCite#1,{%
     \ifx\relax#1%
         \let\next=\relax%
     \else%
         \AddCite{#1}%
         \let\next=\ProcessCite%
     \fi%
     \next%
}

\ndef{\AddCites}[1]{\ProcessCite#1,\relax,}

\ndef{\CiteWithoutExtension}[1]{%
   \AddCites{#1}%
   \LatexCite{#1}%
}

\def\CiteWithExtension[#1]#2{%
   \AddCites{#2}%
   \LatexCite[#1]{#2}%
}

\ndef{\CleverCite}{%
    \ifx\NChar[ %
       \let\MyCite=\CiteWithExtension %
    \else %
       \let\MyCite=\CiteWithoutExtension %
    \fi %
    \MyCite%
}

\renewcommand{\cite}{\futurelet\NChar\CleverCite}
      \ndef{\volume}[1]{{\bf #1}}
      \ndef{\VolYearPP}[3]{\ifnum#2=0 (to appear)\else\volume{#1} (#2), #3\fi}
      \ndef{\VolNoYearPP}[4]{\ifnum#3=0 (to appear)\else\volume{#1} #2 (#3), #4\fi}
      \ndef{\libcode}[1]{}%{{,\bf\ #1}}

\ndef{\jnActaMath}[3]{Acta Math. \VolYearPP{#1}{#2}{#3}}                       % ActM        \libcode{510.5 A18}
\ndef{\jnAdvMath}[3]{Adv.\,in~Math. \VolYearPP{#1}{#2}{#3}}                     % AdvM        \libcode{510.5 A24}
\ndef{\jnAlgAnal}[3]{Algebra i~Analiz \VolYearPP{#1}{#2}{#3}}
\ndef{\jnAmerJMath}[3]{Amer.\,J.\,Math. \VolYearPP{#1}{#2}{#3}}                  % AJM         \libcode{}
\ndef{\jnAmerMathMonth}[3]{Amer.\,Math.\,Monthly \VolYearPP{#1}{#2}{#3}}         % AMM         \libcode{510.5 A3}
\ndef{\jnAnnMath}[4]{Ann. of~Math. \VolNoYearPP{#1}{#2}{#3}{#4}}               % AnnM        \libcode{510.5 A61}
\ndef{\jnAnalMath}[3]{J. Anal. Math. \VolYearPP{#1}{#2}{#3}}                   % JAnalM      \libcode{}
\ndef{\jnArchRatMechAnal}[3]{Arch. Rational Mech. Anal. \VolYearPP{#1}{#2}{#3}}                   % JAnalM      \libcode{}
\ndef{\jnBullLondMathSoc}[3]{Bull. London Math. Soc. \VolYearPP{#1}{#2}{#3}}   % BLMS        \libcode{}
\ndef{\jnBullAMS}[3]{Bull. Amer. Math. Soc. \VolYearPP{#1}{#2}{#3}}   % BLMS        \libcode{}
\ndef{\jnCanMathBull}[3]{Canad. Math. Bull. \VolYearPP{#1}{#2}{#3}}            % CMB         \libcode{}
\ndef{\jnCanMath}[3]{Canad. J.~Math. \VolYearPP{#1}{#2}{#3}}             % CJM         \libcode{}
\ndef{\jnCommMathPhys}[3]{Comm. Math. Phys. \VolYearPP{#1}{#2}{#3}}             % CMP         \libcode{530.1505 C73}
\ndef{\jnCommPDE}[3]{Comm. Partial Differential Equations \VolYearPP{#1}{#2}{#3}}             % CMP         \libcode{530.1505 C73}
\ndef{\jnComptRendue}[3]{C.\,R.~Acad. Sci. Paris S\'er. A-B \VolYearPP{#1}{#2}{#3}}      % CR   \libcode{505 A 161}
\ndef{\jnContMath}[3]{Contemporary Math. \VolYearPP{#1}{#2}{#3}}               %
\ndef{\jnDukeMJ}[3]{Duke Math. J. \VolYearPP{#1}{#2}{#3}}
\ndef{\jnDiffGeom}[3]{J.~Diff. Geom. \VolYearPP{#1}{#2}{#3}}                   % JDG         \libcode{516.3605 J86}
\ndef{\jnErgodicTheory}[3]{Ergodic Theory and Dynamical Systems \VolYearPP{#1}{#2}{#3}} % ETDS  \libcode{}
\ndef{\jnFuncAnal}[3]{J.~Functional Analysis \VolYearPP{#1}{#2}{#3}}           % JFA         \libcode{515.05 J88}
\ndef{\jnFunkAnalPril}[4]{Funct. Anal. Appl. \VolNoYearPP{#1}{#2}{#3}{#4}}  % JFAP
\ndef{\jnGAFA}[3]{GAFA \VolYearPP{#1}{#2}{#3}}                                 % GAFA        \libcode{}
\ndef{\jnIHES}[3]{IHES Publ. Math. (Paris) \VolYearPP{#1}{#2}{#3}}             % IHES        \libcode{}
\ndef{\jnIEOT}[3]{Integral Equations Operator Theory   \VolYearPP{#1}{#2}{#3}} %         \libcode{515.4505 I61}
\ndef{\jnIsrMath}[3]{Israel J.~Math. \VolYearPP{#1}{#2}{#3}}                   % IMJ         \libcode{}
\ndef{\jnKTheory}[3]{K-Theory \VolYearPP{#1}{#2}{#3}}                          % KT          \libcode{}
\ndef{\jnLetMathPhys}[3]{Lett. Math. Phys. \VolYearPP{#1}{#2}{#3}}             % LMP         \libcode{}
\ndef{\jnMathAnn}[3]{Math. Ann. \VolYearPP{#1}{#2}{#3}}                        % MAnn        \libcode{510.5 M51}
\ndef{\jnMathAnalAppl}[3]{J.~Math. Anal. and Appl. \VolYearPP{#1}{#2}{#3}}     % MathAnalAppl  \libcode{510.5 J88} only in repository
\ndef{\jnMathNachr}[3]{Math.\,Nachr. \VolYearPP{#1}{#2}{#3}}
\ndef{\jnMathPhys}[3]{J. Math. Phys. \VolYearPP{#1}{#2}{#3}}
\ndef{\jnMathSocJap}[3]{J. Math. Soc. Japan \VolYearPP{#1}{#2}{#3}}
\ndef{\jnOperTheory}[3]{J.~Operator Theory \VolYearPP{#1}{#2}{#3}}             % JOT         \libcode{}
\ndef{\jnPacJMath}[3]{Pacific J.~Math. \VolYearPP{#1}{#2}{#3}}                  % PJM         \libcode{510.5 P11}
\ndef{\jnPositivity}[3]{Positivity \VolYearPP{#1}{#2}{#3}}
\ndef{\jnProcAmerMS}[3]{Proc. Amer. Math. Soc. \VolYearPP{#1}{#2}{#3}}         % PAMS        \libcode{510.5 A5p}
\ndef{\jnProcCambPhilSoc}[3]{Math. Proc. Camb. Phil. Soc. \VolYearPP{#1}{#2}{#3}}
\ndef{\jnReineAngew}[3]{J.~Reine Angew. Math. \VolYearPP{#1}{#2}{#3}}          % JRAM        \libcode{510.5 J86}
\ndef{\jnTokyoMath}[3]{Tokyo J.~Math. \VolYearPP{#1}{#2}{#3}}
\ndef{\jnTopology}[3]{Topology \VolYearPP{#1}{#2}{#3}}
\ndef{\jnTransAmerMathSoc}[3]{Trans. Amer. Math. Soc. \VolYearPP{#1}{#2}{#3}}
\ndef{\jnIzvANSSSR}[3]{Izv. Akad. Nauk SSSR, Ser. Mat. \VolYearPP{#1}{#2}{#3}}
\ndef{\jnIzvVyshUchZav}[3]{Izv. Vyssh. Uch. Zav., Mat. \VolYearPP{#1}{#2}{#3} (Russian)}
\ndef{\jnIzdatLenUniv}[2]{Izdat. Leningrad. Univ., Leningrad, (#1), #2 (Russian)}
\ndef{\jnFieldsInsComm}[3]{Fields Inst. Comm. \VolYearPP{#1}{#2}{#3}}
\ndef{\jnDoklANSSSR}[3]{Dokl. Akad. Nauk SSSR \VolYearPP{#1}{#2}{#3}}
\ndef{\jnMatZametki}[3]{Matem. zametki \VolYearPP{#1}{#2}{#3}}
\ndef{\jnRussMathSurvey}[3]{Russian Math. Surveys \VolYearPP{#1}{#2}{#3}}
\ndef{\jnSibMathJ}[3]{Sib. Math.~J. \VolYearPP{#1}{#2}{#3}}
\ndef{\jnSovMath}[3]{J.~Soviet math. \VolYearPP{#1}{#2}{#3}}
\ndef{\jnTransMoscMathSoc}[3]{Trans. Moscow Math. Soc. \VolYearPP{#1}{#2}{#3}}
\ndef{\jnUMN}[3]{Uspekhi Mat. Nauk \VolYearPP{#1}{#2}{#3}}

\ndef{\bkTransMathMon}[2]{Trans. Math. Monographs, AMS, \volume{#1}, #2}

\ndef{\pbBirkhauser}[1]{Birkh\"auser, Boston, #1}
\ndef{\pbFactorial}[1]{Moscow, Factorial, #1}
\ndef{\pbGauthier}[1]{Gauthier-Villars, Paris, #1}
\ndef{\pbNauka}[1]{Moscow, Nauka, #1 (Russian)}
\ndef{\pbNaukaR}[1]{Москва, Наука, #1}
\ndef{\pbPrinceton}[1]{Princeton University Press, Princeton, New Jersey, #1}
\ndef{\pbPublPerish}[1]{Publish or Perish Inc., Berkeley, #1}
\ndef{\pbSpringer}[1]{Springer-Verlag, #1}

\ndef{\myauthor}[1]{\mbox{#1}}

\ndef{\Agmon}{\myauthor{Sh.\,Agmon}}
\ndef{\Ahiezer}{\myauthor{N.\,I.\,Ahiezer}}
\ndef{\Arazy}{\myauthor{J.\,Arazy}}
\ndef{\Aronszajn}{\myauthor{N.\,Aronszajn}}
\ndef{\Astashkin}{\myauthor{S.\,V.\,Astashkin}}
\ndef{\Atiyah}{\myauthor{M.\,Atiyah}}
\ndef{\Avron}{\myauthor{J.\,E.\,Avron}}
\ndef{\Azamov}{\myauthor{N.\,A.\,Azamov}}
\ndef{\Banach}{\myauthor{S.\,Banach}}
\ndef{\Benameur}{\myauthor{M-T.\,Benameur}}
\ndef{\Bennett}{\myauthor{C.\,Bennett}}
\ndef{\Berezin}{\myauthor{F.\,A.\,Berezin}}
\ndef{\Berline}{\myauthor{N.\,Berline}}
\ndef{\Birman}{\myauthor{M.\,Sh.\,Birman}}
\ndef{\Blackadar}{\myauthor{B.\,Blackadar}}
\ndef{\Bogolyubov}{\myauthor{N.\,N.\,Bogolyubov}}
\ndef{\Bonsall}{\myauthor{F.\,F.\,Bonsall}}
\ndef{\Bony}{\myauthor{J.\,F.\,Bony}}
\ndef{\BoosBavnbek}{\myauthor{B.\,Boo$\beta$-Bavnbek}}
\ndef{\Bott}{\myauthor{R.\,Bott}}
\ndef{\Branges}{\myauthor{L.\,de Branges}}
\ndef{\Bratteli}{\myauthor{O.\,Bratteli}}
\ndef{\Bredon}{\myauthor{G.\,E.\,Bredon}}
\ndef{\Breuer}{\myauthor{M.\,Breuer}}
\ndef{\Brown}{\myauthor{L.\,G.\,Brown}}
\ndef{\Bruneau}{\myauthor{V.\,Bruneau}}
\ndef{\Buslaev}{\myauthor{V.\,S.\,Buslaev}}
\ndef{\Carey}{\myauthor{A.\,L.\,Carey}}
\ndef{\CareyRW}{\myauthor{R.\,W.\,Carey}} %Richard
\ndef{\Cartan}{\myauthor{H.\,Cartan}}
\ndef{\Chilin}{\myauthor{V.\,I.\,Chilin}}
\ndef{\Coburn}{\myauthor{L.\,A.\,Coburn}}
\ndef{\Connes}{\myauthor{A.\,Connes}}
\ndef{\Cornfeld}{\myauthor{I.\,P.\,Cornfeld}}
\ndef{\Daletskii}{\myauthor{Yu.\,L.\,Daletski\u\i}}   %Daletskii
\ndef{\Dixmier}{\myauthor{J.\,Dixmier}}
\ndef{\DoddsPG}{\myauthor{P.\,G.\,Dodds}}
\ndef{\DoddsTK}{\myauthor{T.\,K.\,Dodds}}
\ndef{\Douglas}{\myauthor{R.\,G.\,Douglas}}
\ndef{\Dubrovin}{\myauthor{B.\,A.\,Dubrovin}}
\ndef{\Dugundji}{\myauthor{J.\,Dugundji}}
\ndef{\Duncan}{\myauthor{J.\,Duncan}}
\ndef{\Dunford}{\myauthor{N.\,Dunford}}
\ndef{\Dykema}{\myauthor{K.\,J.\,Dykema}}
\ndef{\Edwards}{\myauthor{R.\,E.\,Edwards}}
\ndef{\Eilenberg}{\myauthor{S.\,Eilenberg}}
\ndef{\Entina}{\myauthor{S.\,B.\,\`Entina}}
\ndef{\Fack}{\myauthor{T.\,Fack}} %Thierry
\ndef{\Faddeev}{\myauthor{L.\,D.\,Faddeev}}
\ndef{\Farber}{\myauthor{M.\,Farber}}
\ndef{\Farforovskaya}{\myauthor{Yu.\,B.\,Farforovskaya}}
\ndef{\Federer}{\myauthor{H.\,Federer}}
\ndef{\Fedosov}{\myauthor{B.\,V.\,Fedosov}}
\ndef{\Figiel}{\myauthor{T.\,Figiel}} %Tadeush?
\ndef{\Figueroa}{\myauthor{H.\,Figueroa}}
\ndef{\Fillmore}{\myauthor{P.\,A.\,Fillmore}}
\ndef{\Fomenko}{\myauthor{A.\,T.\,Fomenko}} % Anatolii Trofimovich
\ndef{\Fomin}{\myauthor{S.\,V.\,Fomin}}
\ndef{\Frohlich}{\myauthor{J.\,Fr\"ohlich}}
\ndef{\Fuglede}{\myauthor{B.\,Fuglede}}
\ndef{\Furutani}{\myauthor{K.\,Furutani}}
\ndef{\Gelfand}{\myauthor{I.\,M.\,Gelfand}}
\ndef{\Gesztesy}{\myauthor{F.\,Gesztesy}}     %Fritz
\ndef{\Getzler}{\myauthor{E.\,Getzler}} % Ezra
\ndef{\Gilkey}{\myauthor{P.\,B.\,Gilkey}}
\ndef{\Gitler}{\myauthor{S.\,Gitler}}
\ndef{\Glazman}{\myauthor{I.\,M.\,Glazman}}
\ndef{\Glimm}{\myauthor{J.\,Glimm}}
\ndef{\Gohberg}{\myauthor{I.\,C.\,Gohberg}}
\ndef{\Goldshtein}{\myauthor{Ya.\,Goldshtein}}
\ndef{\Golze}{\myauthor{F.\,Golze}}
\ndef{\GraciaBondia}{\myauthor{J.\,M.\,Gracia-Bond\'{i}a}}
\ndef{\Greenleaf}{\myauthor{F.\,P.\,Greenleaf}}
\ndef{\Gromov}{\myauthor{M.\,Gromov}}
\ndef{\Gunning}{\myauthor{R.\,C.\,Gunning}}
\ndef{\Haagerup}{\myauthor{U.\,Haagerup}}
\ndef{\Haag}{\myauthor{R.\,Haag}}
\ndef{\Halmos}{\myauthor{P.\,R.\,Halmos}}
\ndef{\Hardy}{\myauthor{G.\,H.\,Hardy}}
\ndef{\Herbst}{\myauthor{I.\,W.\,Herbst}}
\ndef{\Higson}{\myauthor{N.\,Higson}}  % Nigel
\ndef{\Hoermander}{\myauthor{L.\,H\"ormander}} % Lars
\ndef{\Hoffman}{\myauthor{K.\,Hoffman}} % Kenneth Hoffman
\ndef{\Ito}{\myauthor{K.\,Ito}}
\ndef{\Ikebe}{\myauthor{T.\,Ikebe}}
\ndef{\Jaffe}{\myauthor{A.\,Jaffe}}
\ndef{\James}{\myauthor{I.\,M.\,James}}
\ndef{\Javrjan}{\myauthor{V.\,A.\,Javrjan}}
\ndef{\Jitomirskaya}{\myauthor{S.\,Jitomirskaya}}
\ndef{\Kadison}{\myauthor{R.\,V.\,Kadison}}
\ndef{\Kalton}{\myauthor{N.\,J.\,Kalton}} % Nigel
\ndef{\Kato}{\myauthor{T.\,Kato}} % Tosio
\ndef{\Kobayashi}{\myauthor{S.\,Kobayashi}}
\ndef{\Koplienko}{\myauthor{L.\,S.\,Koplienko}}
\ndef{\Korotyaev}{\myauthor{E.\,Korotyaev}}
\ndef{\Kosaki}{\myauthor{H.\,Kosaki}}
\ndef{\Kostrykin}{\myauthor{V.\,Kostrykin}}
\ndef{\Kotani}{\myauthor{S.\,Kotani}}
\ndef{\Krein}{\myauthor{Kre\u\i n}}
\ndef{\KreinMG}{\myauthor{M.\,G.\,Kre\u\i n}}
\ndef{\KreinSG}{\myauthor{S.\,G.\,Kre\u\i n}}
\ndef{\Kuroda}{\myauthor{S.\,T.\,Kuroda}}
\ndef{\Leichtnam}{\myauthor{E.\,Leichtnam}}
\ndef{\Lesch}{\myauthor{M.\,Lesch}}
\ndef{\Lesniewski}{\myauthor{A.\,Lesniewski}}
\ndef{\Levitan}{\myauthor{B.\,M.\,Levitan}}
\ndef{\Lidskii}{\myauthor{V.\,B.\,Lidskii}}
\ndef{\Lifshitz}{\myauthor{I.\,M.\,Lifshitz}}
\ndef{\Lindenstrauss}{\myauthor{J.\,Lindenstrauss}}
\ndef{\Loday}{\myauthor{J.-L.\,Loday}}
\ndef{\Lord}{\myauthor{S.\,Lord}}      %Steven
\ndef{\Lorentz}{\myauthor{G.\,Lorentz}}
\ndef{\Magnus}{\myauthor{W.\,Magnus}}
\ndef{\Makarov}{\myauthor{K.\,A.\,Makarov}}
\ndef{\MakarovN}{\myauthor{N.\,Makarov}}
\ndef{\Mathai}{\myauthor{V.\,Mathai}}         %Varghese?
\ndef{\McKean}{\myauthor{H.\,P.\,McKean}}
\ndef{\Mishchenko}{\myauthor{A.\,S.\,Mishchenko}}
\ndef{\Molchanov}{\myauthor{S.\,A.\,Molchanov}}
\ndef{\Moore}{\myauthor{C.\,C.\,Moore}}
\ndef{\Moscovici}{\myauthor{H.\,Moscovici}}  %Henri?
\ndef{\Motovilov}{\myauthor{A.\,K.\,Motovilov}}
\ndef{\Moyer}{\myauthor{R.\,D.\,Moyer}}
\ndef{\Naboko}{\myauthor{S.\,N.\,Naboko}}
\ndef{\Narasimhan}{\myauthor{R.\,Narasimhan}}
\ndef{\Nomizu}{\myauthor{K.\,Nomizu}}
\ndef{\Novikov}{\myauthor{S.\,P.\,Novikov}}
\ndef{\Osterwalder}{\myauthor{K.\,Osterwalder}}
\ndef{\Patodi}{\myauthor{V.\,Patodi}}
\ndef{\Pagter}{\myauthor{B.\,de~Pagter}}  %Ben
\ndef{\Pastur}{\myauthor{L.\,A.\,Pastur}}  %Ben
\ndef{\Pavlov}{\myauthor{B.\,S.\,Pavlov}}
\ndef{\Pedersen}{\myauthor{G.\,K.\,Pedersen}}
\ndef{\Peller}{\myauthor{V.\,V.\,Peller}}
\ndef{\Perera}{\myauthor{V.\,S.\,Perera}}
\ndef{\Petunin}{\myauthor{Ju.\,I.\,Petunin}}
\ndef{\Phillips}{\myauthor{J.\,Phillips}}  %John
\ndef{\Piazza}{\myauthor{P.\,Piazza}}   %Paolo
\ndef{\Pincus}{\myauthor{J.\,D.\,Pincus}}   %Joel
\ndef{\Poincare}{Poincar\'e}
\ndef{\Postnikov}{\myauthor{M.\,M.\,Postnikov}} % Mikhail
\ndef{\Povzner}{\myauthor{A.\,Ya.\,Povzner}}
\ndef{\Prinzis}{\myauthor{R.\,Prinzis}}
\ndef{\Privalov}{\myauthor{I.\,I.\,Privalov}}
\ndef{\Pushnitski}{\myauthor{A.\,B.\,Pushnitski}} % Alexander
\ndef{\Raeburn}{\myauthor{I.\,Raeburn}}
\ndef{\Raikov}{\myauthor{G.\,Raikov}}
\ndef{\Reed}{\myauthor{M.\,Reed}}
\ndef{\Rennie}{\myauthor{A.\,Rennie}}
\ndef{\Rickart}{\myauthor{C.\,E.\,Rickart}}
\ndef{\Riesz}{\myauthor{F.\,Riesz}}
\ndef{\Ringrose}{\myauthor{J.\,Ringrose}}
\ndef{\Rio}{\myauthor{R.\,del Rio}}
\ndef{\Robinson}{\myauthor{D.\,Robinson}}
\ndef{\Rossi}{\myauthor{H.\,Rossi}}
\ndef{\Rudin}{\myauthor{W.\,Rudin}}
\ndef{\Ruelle}{\myauthor{D.\,Ruelle}}
\ndef{\Ruzhansky}{\myauthor{M.\,Ruzhansky}}
\ndef{\Sakai}{\myauthor{Sh.\,Sakai}}
\ndef{\Sargsjan}{\myauthor{I.\,S.\,Sargsjan}}
\ndef{\Sato}{\myauthor{H.\,Sato}}
\ndef{\Schaeffer}{\myauthor{D.\,G.\,Schaeffer}}
\ndef{\Schluchtermann}{\myauthor{G.\,Schluchtermann}}
\ndef{\Schochet}{\myauthor{C.\,Schochet}}
\ndef{\SchroedingerE}{\myauthor{E.\,Schr\"odinger}}
\ndef{\Schroedinger}{\myauthor{Schr\"odinger}}
\ndef{\Schrohe}{\myauthor{E.\,Schrohe}}
\ndef{\Schwartz}{\myauthor{J.\,T.\,Schwartz}}
\ndef{\Sedaev}{\myauthor{A.\,A.\,Sedaev}}
\ndef{\Seiler}{\myauthor{R.\,Seiler}}
\ndef{\Semenov}{\myauthor{E.\,M.\,Semenov}}
\ndef{\Shabat}{\myauthor{B.\,V.\,Shabat}}
\ndef{\Shafarevich}{\myauthor{I.\,R.\,Shafarevich}}
\ndef{\Sharpley}{\myauthor{R.\,Sharpley}}
\ndef{\Shilov}{\myauthor{G.\,E.\,Shilov}}
\ndef{\Shirkov}{\myauthor{D.\,V.\,Shirkov}}
\ndef{\Shubin}{\myauthor{M.\,A.\,Shubin}}
\ndef{\Silverman}{\myauthor{H.\,Silverman}}
\ndef{\Simon}{\myauthor{B.\,Simon}}
\ndef{\Sinai}{\myauthor{Ya.\,G.\,Sinai}}
\ndef{\Singer}{\myauthor{I.\,M.\,Singer}}
\ndef{\Solomyak}{\myauthor{M.\,Z.\,Solomyak}}
\ndef{\Soloviev}{\myauthor{Yu.\,P.\,Soloviev}}
\ndef{\Spivak}{\myauthor{M.\,Spivak}}
\ndef{\Stein}{\myauthor{E.\,M.\,Stein}}
\ndef{\Stenkin}{\myauthor{V.\,V.\,Sten'kin}}
\ndef{\Stratila}{\myauthor{S.\,Stratila}}
\ndef{\Sucheston}{\myauthor{L.\,Sucheston}}
\ndef{\Sukochev}{\myauthor{F.\,A.\,Sukochev}}
\ndef{\Switzer}{\myauthor{R.\,M.\,Switzer}}
\ndef{\SzNagy}{\myauthor{B.\,Sz.-Nagy}}
\ndef{\Takesaki}{\myauthor{M.\,Takesaki}}
\ndef{\Taylor}{\myauthor{M.\,E.\,Taylor}}
\ndef{\Treves}{\myauthor{F.\,Treves}}
\ndef{\Troitsky}{\myauthor{E.\,V.\,Troitsky}}
\ndef{\Tzafriri}{\myauthor{L.\,Tzafriri}}
\ndef{\Varilly}{\myauthor{J.\,C.\,V\'{a}rilly}}
\ndef{\Vergne}{\myauthor{M.\,Vergne}}
\ndef{\Vladimirov}{\myauthor{V.\,S.\,Vladimirov}}
\ndef{\Voiculescu}{\myauthor{D.\,Voiculescu}}
\ndef{\Weiss}{\myauthor{G.\,Weiss}}
\ndef{\Wells}{\myauthor{R.\,O.\,Wells}}
\ndef{\Williams}{\myauthor{J.\,P.\,Williams}}
\ndef{\Winkler}{\myauthor{S.\,Winkler}}
\ndef{\Witten}{\myauthor{E.\,Witten}}
\ndef{\Wodzicki}{\myauthor{M.\,Wodzicki}}
\ndef{\Wojciechowski}{\myauthor{K.\,P.\,Wojciechowski}}
\ndef{\Yafaev}{\myauthor{D.\,R.\,Yafaev}}
\ndef{\Yosida}{\myauthor{K.\,Yosida}}
\ndef{\Zsido}{\myauthor{L.\,Zsido}}

%% file: MyListOfRef.tex
% MyListOfRef.tex

\mathsurround 0pt

% End of MyListOfRef.tex